\documentclass{article}
\usepackage{amsmath, amsfonts, amssymb}

\newtheorem{theorem}{Theorem}
\newtheorem{lemma}[theorem]{Lemma}

\newcommand{\convex}{\mathrm{conv}\,}
\newcommand{\volume}{\mathrm{vol}}

\newcommand{\vol}[1]{\volume(#1)}
\newcommand{\R}{\mathbb{R}}
\newcommand{\ipr}[2]{\langle #1, #2 \rangle}
\newcommand{\norm}[1]{\lVert#1\rVert}
\newcommand{\abs}[1]{\lvert#1\rvert}

\providecommand{\qed}{\hfill \rule{1ex}{1ex}}

\begin{document}

\title{Sets of unit vectors with small pairwise sums}
\author{Konrad J.\ Swanepoel \thanks{I thank the referee for remarks leading to an improved paper.}\\
Department of Mathematics and Applied Mathematics \\ 
University of Pretoria \\ Pretoria 0002 \\ South Africa \\ 
e-mail: \texttt{konrad@math.up.ac.za}}
\date{}
\maketitle

\begin{abstract}
We study the sizes of \emph{$\delta$-additive} sets of unit vectors in a $d$-dimensional normed space: the sum of any two vectors has norm at most $\delta$.
One-additive sets originate in finding upper bounds of vertex degrees of Steiner Minimum Trees in finite dimensional smooth normed spaces (Z.~F{\"u}redi, J.~C. Lagarias, F.~Morgan, 1991).
We show that the maximum size of a $\delta$-additive set over all normed spaces of dimension $d$ grows exponentially in $d$ for fixed $\delta>2/3$, stays bounded for $\delta<2/3$, and grows linearly at the threshold $\delta=2/3$.
Furthermore, the maximum size of a $2/3$-additive set in $d$-dimensional normed space has the sharp upper bound of $d$, with the single exception of spaces isometric to three-dimensional $\ell^1$ space, where there exists a $2/3$-additive set of four unit vectors.

\end{abstract}

\section{Introduction}\label{intro}
Let $X$ be a real normed space of dimension $d\geq 1$, with norm $\norm{\cdot}$.
Let $0<\delta< 2$.
A {\em $\delta$-additive set} in $X$ is a set of unit vectors $S$ satisfying $\norm{x+y}\leq\delta$ for distinct $x,y\in S$.
We define $N_X(\delta)$ to be the largest cardinality of a $\delta$-additive set in $X$.
These notions originate in the analysis of geometric Steiner Minimum Trees in combinatorial optimization.
If $X$ is smooth, then $N_{X^\ast}(1)$ is an upper bound for the maximum degree in any Steiner Minimal Tree in $X$, where $X^\ast$ is the dual of $X$; see \cite{LM, S2}.
It is easily seen that $N_X(\delta)$ is finite for all $\delta\in(0,2)$.
In this note we investigate the maximum of $N_X(\delta)$ over all $X$ of a fixed dimension, keeping $\delta$ fixed.
Let $N_d(\delta)=\max N_X(\delta)$, where the maximum is over all $X$ of dimension $d$.

First of all, note that if $0<\delta<2/3$, then the triangle inequality immediately gives $N_X(\delta) =2$ for all $X$.
If $\delta>2/3$, the following generalization of \cite[Theorem 2.4]{FLM} shows that $N_d(\delta)$ grows at least exponentially in $d$.

\begin{theorem}\label{normexp}
For all $\delta>\tfrac{2}{3}$ there exist $\epsilon=\epsilon(\delta)>0$ such that for all sufficiently large $d$ there exists a $d$-dimensional normed space $X$ such that $N_X(\delta) > (1+\epsilon)^d$.
\end{theorem}

By packing arguments it is easily seen that there is also an upper bound for $N_d(\delta)$ exponential in $d$ for fixed $\delta$.
The following is the best such an upper bound we have.

\begin{theorem}\label{upperbound}
If $X$ is a $d$-dimensional normed space, then
\[
N_{X}(\delta) \leq 2\left(\frac{2}{2-\delta}\right)^d.
\]
\end{theorem}

The proof (in Section~\ref{upperboundproof}) uses the Brunn-Minkowski inequality.

Surprisingly, in the remaining case of $\delta=2/3$, $N_d(\delta)$ grows linearly in $d$.
The following theorem gives the exact value of $N_d(\delta)$, and shows that for $d=3$, spaces isometric too three-dimensional $\ell_1$ are exceptional.
This theorem can therefore be considered to be a characterization of the three-dimensional affine regular octahedron in the collection of all centrally symmetric convex bodies of any finite dimension.

\begin{theorem}\label{linear}
Let $X$ be a $d$-dimensional normed space.

If $d\neq 1, 3$, then $N_X(2/3) \leq d$, equality being attained e.g.\ if the unit ball of $X$ is a cube.

If $d=3$, then $N_X(2/3) \leq 4$, with equality iff the unit ball of $X$ is an affine regular octahedron.
\end{theorem}

In the sequel we fix our $d$-dimensional space to be $\R^d$ with inner product $\ipr{\cdot}{\cdot}$.
For $S\subseteq\R^d$ let $\convex S$ be the convex hull of $S$.
If $S$ is Lebesgue measurable, let $\vol{S}$ be the Lebesgue measure or \emph{volume} of $S$.
Let $\ell_\infty^d$ be $\R^d$ with norm $\norm{(x_1,\dots,x_d)}_\infty=\max_i\abs{x_i}$, and $\ell_1^d$ be $\R^d$ with norm $\norm{(x_1,\dots,x_d)}_1=\abs{x_1}+\dots+\abs{x_d}$.
Note that the unit ball of $\ell_\infty^d$ is a $d$-dimensional cube.
Also note that a three-dimensional normed space is isometric to $\ell_3^1$ iff its unit ball is an \emph{affine regular octahedron} centered at $0$, i.e.\ a non-singular linear image of the regular octahedron centered at $0$.
 
We now state a technical lemma important to the proofs of Theorems~\ref{normexp} and \ref{linear} (sections~\ref{normexpproof} and \ref{linearproof}, respectively).
The lemma reduces the existence of norms admitting a given $\delta$-additive set to a set of linear inequalities.
The proof is in Section~\ref{lemmaproof}.

\begin{lemma}\label{normexist}
Let $x_1,\dots,x_m\in\R^d\setminus\{0\}$ and $0<\delta\leq 2$.
There exists a norm $\norm{\cdot}$ on $\R^d$ such that
\begin{equation}\label{codeexist}
\begin{split}
\norm{x_i} = 1 & \;\forall\, i=1,\dots,m,\\
\norm{x_i+x_j}\leq\delta & \;\forall\, 1\leq i<j\leq m,
\end{split}
\end{equation}
iff there exist $y_1,\dots,y_m\in\R^d$ such that
\begin{equation}\label{yexist}
\begin{split}
\ipr{x_i}{y_i} = 1 & \;\forall\, i=1,\dots,m,\\
-1\leq \ipr{x_j}{y_i}\leq\delta-1&\;\forall\, \text{distinct } i,j=1,\dots,m,\\
-\delta\leq\ipr{x_j+x_k}{y_i}&\;\forall\, i,j,k=1,\dots,m\text{ with } j\neq k.
\end{split}
\end{equation}
\end{lemma}

\section{Proof of Theorem~\ref{upperbound}}\label{upperboundproof}
Let $S$ be a $\delta$-additive set containing $N$ vectors.
Then for distinct $x,y\in S$, $\norm{x+y}\leq\delta, \norm{x-y}\geq 2-\delta$.
We partition $S$ into two sets $S_1$ and $S_2$ of sizes $\lfloor N/2\rfloor$ and $\lceil N/2\rceil$, respectively.
Let $V_i := \bigcup_{x\in S_i} B(x,1-\delta/2)$ for $i=1,2$.
Each $V_i$ consists of closed balls with disjoint interiors, and therefore, $\vol{V_1}=(\lfloor N/2\rfloor)2^{-d}\vol{B}$ and $\vol{V_2}=(\lceil N/2\rceil)2^{-d}\vol{B}$.
Also, $V_1+V_2\subseteq B(0,3-\delta)$.
By the Brunn-Minkowski inequality (see \cite{BZ}) we obtain
\[(\lfloor N/2\rfloor^{1/d} + \lceil N/2\rceil^{1/d})(1-\delta/2)\leq 3-\delta,\]
and the result follows.
\qed

\section{Proof of Theorem~\ref{linear}}\label{linearproof}
To see that equality may hold for $d\geq 4$, consider the space $X=\ell_\infty^d$ with unit ball $[-1,1]^d$, and let $S$ be the set of all coordinate permutations of 
\[(1,-\tfrac{1}{3}, -\tfrac{1}{3},\dots, -\tfrac{1}{3})\in\ell_\infty^d.\]
Then $S$ is a set of $d$ unit vectors and $\norm{x+y}=\tfrac{2}{3}$ for distinct $x,y\in S$.

For $d=3$, note that the 4 vectors of $S$ in $\ell_\infty^4$ are all in the hyperplane $\{x\in\R^4:\sum_i x_i=0\}$, and thus in a 3-dimensional subspace.
It is easy to see that this subspace is isometric to $\ell_1^3$, with unit ball the octahedron \[\convex\{(\pm 1,0,0), (0,\pm 1,0), (0,0,\pm 1)\}.\]

We first derive the upper bound for $d\neq 1,3$.
Let $\{x_1,\dots,x_m\}$ be a $\frac{2}{3}$-additive set of unit vectors in $X$.
We suppose for the sake of contradiction that $m=d+1$.
By Lemma~\ref{normexist} there are $y_1,\dots,y_m\in X$ such that
$$\ipr{x_j}{y_i}\leq-\tfrac{1}{3}\leq\tfrac{1}{2}\ipr{x_j+x_k}{y_i}\;\forall\, i,j,k \text{ with } j\neq k,j\neq i.$$
Thus $\ipr{x_j}{y_i}\leq \ipr{x_k}{y_i}$.
But similarly, $\ipr{x_k}{y_i}\leq\ipr{x_j}{y_i}\;\forall\, i,j,k$ with $k\neq j, k\neq i$.
Thus  $\ipr{x_j}{y_i} = \ipr{x_k}{y_i} = -\tfrac{1}{3} \;\forall \text{ distinct } i,j,k$, i.e.\
\begin{equation}\label{dotpr}
\ipr{x_j}{y_i} = \begin{cases} 1 & \text{ if }i=j,\\ -\tfrac{1}{3}& \text{ if }i\neq j.\end{cases}
\end{equation}
Note that we have used the fact $m\geq 3$.

Suppose that $\lambda_i$ is a sequence of scalars for which $\sum_{i=1}^m\lambda_i x_i = 0$.
Then for all $j$ we have
$$0=\sum_{i=1}^m\lambda_i\ipr{x_i}{y_j} = \tfrac{4}{3}\lambda_j - \tfrac{1}{3}\sum_{i=1}^m\lambda_i.$$
Thus all $\lambda_j$'s are equal: 
\begin{equation}\label{linind}
\lambda_j=:\lambda\;\forall j=1,\dots, m, 
\end{equation}
and $(\tfrac{4}{3}-\tfrac{1}{3}m)\lambda = 0$.
Thus $x_1,\dots,x_m$ are linearly independent, since $m\neq 4$, contradicting $m=d+1$.

We now treat the remaining case $d=3$ (as $d=1$ is trivial). 
By considering $X$ as a subspace of some 4-dimensional space, we immediately obtain from the previous argument that $N_X(2/3)\leq 4$.
Alternatively, we can argue directly as follows.
Using the John-Loewner ellipsoid (see e.g.\ \cite[Theorem 3.3.6]{Thompson}), we obtain an inner product $\ipr{\cdot}{\cdot}$ and norm $\norm{\cdot}_2=\sqrt{\ipr{\cdot}{\cdot}}$ in $X$ such that $$\norm{x}\leq\norm{x}_2\leq\sqrt{3}\norm{x}.$$
Let $S=\{x_1,\dots,x_m\}$ be a $\frac{2}{3}$-additive set of unit vectors in $X$.
Then $\norm{x}_2\geq 1$ for all $x\in S$, and $\norm{x+y}_2^2\leq 3\norm{x+y}^2\leq \frac{4}{3}$ for all distinct $x,y\in S$.
Thus $\tfrac{4}{3}\geq 2+2\ipr{x}{y}$, and $\ipr{x}{y}\leq -\frac{1}{3}$.
It follows that $\norm{\sum_{i=1}^m x_i}_2^2 \leq m-\frac{2}{3}\binom{m}{2}$, implying $m\leq 4$.

It remains to show that $\ell_1^3$ is the unique $3$-dimensional Minkowski space admitting four unit vectors satisfying $\norm{x_i+x_j}\leq \frac{2}{3}$ for $i\neq j$.
Suppose we have four such vectors for some norm $\norm{\cdot}$.
From the triangle inequality it follows that $\norm{x_i+x_j}=\frac{2}{3}$ for $i\neq j$, i.e.\ $z_i := \frac{3}{2}(x_i+x_4)$, $i=1,2,3$ are unit vectors that are furthermore linearly independent:
If $\sum_i\lambda_i z_i=0$ then $\sum_i\lambda_i\ipr{z_i}{y_j}=0$ for each $j$, and from \eqref{dotpr} it follows after some calculation that all $\lambda_i=0$.

Since $x_1,\dots,x_4$ are linearly dependent, we have $\sum_i\lambda_i x_i=0$ for some $\lambda_i$ not all $0$.
From \eqref{linind} it follows that the $\lambda_i$ are equal, and therefore,
$\sum_i x_i = 0$.
Hence $z_1+z_2+z_3 = \frac{3}{2}(x_1+x_2+x_3+3x_4)=3x_4$, and $\norm{\frac{1}{3}(z_1+z_2+z_3)}=1$.
Thus, $\convex\{z_1,z_2,z_3\}$ is a face of the unit ball.

Similarly, $\norm{\frac{1}{3}(\epsilon_1 z_1+\epsilon_2 z_2 +\epsilon_3 z_3)}=1$ for any $\epsilon_i=\pm 1$, and 
$$\convex\{\epsilon_1 z_1,\epsilon_2 z_2,\epsilon_3 z_3\}, \epsilon_i=\pm 1$$ 
are all faces of the unit ball, which must therefore be the affine regular octahedron $\convex\{\pm z_1,\pm z_2,\pm z_3\}$.
It follows that $\norm{\sum_i\lambda_i z_i} = \norm{(\lambda_1,\lambda_2,\lambda_3)}_1$.

Note that it is easy to calculate
\begin{align*}
x_1 &=\frac{1}{3}(z_1-z_2-z_3), & x_2 &=\frac{1}{3}(-z_1+z_2-z_3),\\
x_3 &=\frac{1}{3}(-z_1-z_2+z_3), & x_4 &=\frac{1}{3}(z_1+z_2+z_3).
\end{align*}
Except for reflections in the coordinate planes or permutations of the coordinates (i.e.\ linear isometries of $\ell_1^3$), these vectors form the unique $2/3$-additive set of four points in $\ell_1^3$.
They are the centroids of four alternate faces of the octahedron, and the vertex set of a regular tetrahedron.
\qed

\section{Proof of Theorem~\ref{normexp}}\label{normexpproof}
Let $\delta':= (3\delta-2)/(6-\delta) > 0$.
(Recall that we assume $\delta < 2$.)
By a result of \cite{Wyner} it is possible to find at least $(1+\epsilon)^d$ euclidean unit vectors $v_i$ such that $\abs{\ipr{v_i}{v_j}}<\delta'$ for all distinct $i,j$.
Regard $\R^{d+1}$ as  $\R^d\oplus\R$, identify $\R^d$ with the $\R^d$ factor, and let $e$ be a unit vector orthogonal to $\R^d$.
Let $x_i:=v_i + e$ for all $i$.
We now check that Lemma~\ref{normexist} may be applied with $y_i:= \lambda v_i + (1-\lambda)e$, where $\lambda = \tfrac{2}{3}-\frac{\delta}{4}$.

Firstly, for all $i$,
$$\ipr{x_i}{y_i}=\lambda+1-\lambda=1.$$

Secondly, for all $i\neq j$,
$$\ipr{x_j}{y_i} = \lambda \ipr{v_i}{v_j}+1-\lambda
	\begin{cases}
	\leq \lambda\delta'+1-\lambda = \delta-1,\\
	\geq-\lambda\delta'+1-\lambda =-\delta/2\geq -1.
	\end{cases}$$
Thus for all distinct $i,j$,
$$\ipr{x_i+x_j}{y_i} = 1+\ipr{x_j}{y_i}\geq 1-1 \geq 0 \geq -\delta.$$

Thirdly, for all $i,j,k$ such that $j\neq i, k\neq i$,
$$\ipr{x_i+x_k}{y_i}\geq-\delta/2-\delta/2=-\delta.$$

Lemma~\ref{normexist} now gives the required norm.
\qed

\section{Proof of Lemma~\ref{normexist}}\label{lemmaproof}
$\Rightarrow$ We choose norming functionals $y_i$ for $x_i$:
$$\ipr{x_i}{y_i} = 1 = \sup_{\norm{x}=1} \ipr{x}{y_i}.$$
Then for $i\neq j$ we have
$$1+\ipr{x_j}{y_i} = \ipr{x_i+x_j}{y_i}\leq\norm{x_i+x_j}\leq\delta,$$
and hence $\ipr{x_j}{y_i}\leq\delta-1$.
Also, $-\ipr{x_j}{y_i}\leq\norm{-x_j}=1$, and therefore, $\ipr{x_j}{y_i}\geq -1$.
Also, for $j\neq k$ and $i$, 
$$\ipr{-x_j-x_k}{y_i}\leq\norm{-x_j-x_k}\leq\delta,$$
and $\ipr{x_j+x_k}{y_i}\geq-\delta.$

$\Leftarrow$ Let $K$ be the convex hull of 
$$\{\pm x_i:i=1,\dots,m\}\cup\{\pm\tfrac{1}{\delta}(x_i+x_j):1\leq i<j\leq m\}.$$
Then $K$ is centrally symmetric and convex.
If $K$ has no interior points (i.e.\ if $K$ is contained is some hyperplane), we ``thicken'' $K$:
Let $V$ be the linear span of $K$, $U=V^\perp$ and $e_1,\dots,e_k$ an orthonormal basis for $U$.
Replace $K$ by $K+[-1,1]e_1+\dots+[-1,1]e_k$.
We now have that $K$ is a centrally symmetric convex body defining a norm
$$\norm{x}_K:= \inf\{t\geq 0: x\in tK\}.$$

Obviously, $\norm{x_i}_K\leq 1 \;\forall i$ and $\norm{x_i+x_j}_K\leq\delta \;\forall\, i\neq j$.
It remains to show that $\norm{x_i}=1$, i.e.\ that each $x_i$ is on the boundary of $K$.
It is sufficient to show that the closed half space $\{x\in\R^d:\ipr{x}{y_i}\leq 1\}$ contains $K$.
By the definition of $K$ it is sufficient to show that
\begin{align}
\pm\ipr{x_j}{y_i}\leq 1 &\;\forall\, i,j\label{eq1}\\
\text{ and } \pm\ipr{x_j+x_k}{y_i}\leq\delta& \;\forall\, i,j,k \text{ with } j\neq k.\label{eq2}
\end{align}
Now \eqref{eq1} holds, since $\ipr{x_i}{y_i}=1$, and for $j\neq i$, $-1\leq\ipr{x_j}{y_i}\leq\delta-1\leq 1.$
Also, $$\ipr{x_i+x_j}{y_i} = 1+\ipr{x_j}{y_i} = 1+\ipr{x_j}{y_i}
	\begin{cases}
	\leq 1+\delta-1=\delta\\
	\geq 1-1>-\delta
	\end{cases}$$
for $i\neq j$, and
$$\ipr{x_j+x_k}{y_i}
	\begin{cases}
	\leq \delta-1+\delta-1\leq\delta\\
	\geq -\delta
	\end{cases}$$
for $j\neq i,k\neq i$,
and so \eqref{eq2} holds.
\qed

\end{document}